\documentclass[12pt,reqno]{article}
\usepackage{graphicx}
\usepackage{amsmath}
\usepackage{amsthm}
\usepackage{latexsym,bm}
\usepackage{amssymb}
\usepackage{indentfirst}
\newtheorem{thm}{Theorem}[section]
\newtheorem{cor}[thm]{Corollary}

\numberwithin{equation}{section}

\usepackage{graphicx}

\begin{document}

\title{{\bf  Good shadows, dynamics, and convex hulls}}

\author{{\bf F. Fontenele\ and  F. Xavier}}

\date{}
\maketitle





\begin{quote}
\small {\bf Abstract}. The Ekeland variational principle implies
what can be regarded as a strong version, in the $C^1$ category,
of the Yau minimum principle: under the appropriate hypotheses
{\it every} minimizing sequence admits a {\it good shadow}, a
second minimizing sequence that has good properties and is
asymptotic to the original one. Using arguments from dynamical
systems, we give another proof of this result and also establish,
with the aid of Gromov's theorem on monotonicity of volume ratios,
a special case of a conjecture claiming the existence of good
shadows in the original $C^2$ setting of the Yau minimum
principle. The interest in having an abundance of good shadows
stems from the fact that this is a desirable property if one wants
to refine the applications of the asymptotic minimum principle, as
it allows for information to be localized at infinity. These ideas
are applied in this paper to the study of the convex hulls of
complete submanifolds of Euclidean $n$-space that have controlled
Grassmanian-valued Gauss maps.

\end{quote}

\section{Introduction}

The well known Yau minimum principle (\cite{CY},\cite
{O},\cite{Y1}), stated below, is a powerful tool in geometric
analysis (see, for instance, \cite{JK},\cite{RRS},\cite{Y2}):

\begin{thm}\label{OY} Let $M$ be a complete manifold
whose Ricci curvature is bounded from below, and $f:M\to\mathbb R$
of class $C^2$ such that $\inf _M f>-\infty$. Then,  there exists
a sequence $(x_n)$ in $M$ satisfying $f(x_n)\to\inf_Mf$, $||\nabla
f(x_n)||\to 0$ and  $\liminf_{n\to\infty}\Delta f(x_n)\geq 0$.
\end{thm}
By analogy, one may ask the following question: If $f$ is supposed
to be only of class $C^1$, is there a minimizing sequence for $f$
along which the gradient of $f$ is small? In the following theorem
we observe not only that such special sequences exist, but also
that they can be found asymptotically near \it any \rm specified
minimizing sequence.

\begin{thm}\label{SS} Let $M$ be a complete manifold, and $f:M\to
\Bbb R$ of class $C^1$ such that $\inf _M f>-\infty$. Then, for
every sequence $(x_n)$ in $M$ such that $f(x_n)\to \inf _M f$,
there exists a sequence $(y_n)$ in $M$ such that $d(x_n,y_n)\to
0$, $f(y_n)\to\inf_Mf$ and $||\nabla f(y_n)||\to 0$.
\end{thm}

As we will see, this abundance of  \lq\lq good\rq\rq  minimizing
sequences will allow us to establish a sharp geometric application
that does not follow from the Yau minimum principle, even in the
$C^{\infty}$ case (Theorem \ref{ConvHull}).

The formal similarity between Theorems \ref{OY} and \ref{SS} is
manifest, except for the statement about special minimizing
sequences that are asymptotically close to an arbitrary minimizing
one.

\vskip10pt

\noindent{\bf Definition.} Let $M$ be a complete non-compact
Riemannian manifold, and $f:M\to \mathbb R$ a function of class
$C^2$ satisfying $\inf f>-\infty$. A minimizing sequence $(p_n)$
of $f$ is said to admit \it a good shadow \rm if there exists  a
minimizing sequence $(q_n)$ in $M$ such that $ \lim d(p_n,
q_n)=0$, \;$\lim ||\nabla f(q_n)|| = 0$, \ and\; $ \liminf \Delta
f(q_n)\geq 0$ (likewise, if $f$ is $C^1$ one only requires the
first two properties).

\vskip10pt

The analogy between the two theorems above will be complete if the
following can be shown to be  true:

\vskip10pt

\noindent \bf {Conjecture}.\rm  \  Let  $f:M\to \mathbb R$ be of
class $C^2$, $\inf f>-\infty$, $M$  a complete manifold with Ricci
curvature  bounded below. Then every minimizing sequence of $f$
has a good  shadow.\rm

\vskip10pt

In $\S2$ we offer two proofs of  Theorem \ref{SS}. The first one
is based on the Ekeland variational principle (\cite {E},
\cite{Fi}, \cite{St}), a well known result in control theory and
non-linear analysis. A new line of argument, based on ideas from
dynamical systems, is also presented. The reason for including the
second proof here is that a considerable elaboration of it,
together with a result of Gromov, yields the following special
case of the good shadows conjecture (see also Theorem
\ref{oscil}):

\begin{thm}\label{weakyau} Let $M^m$ be a complete manifold with Ricci curvature bounded from below,
and $f:M\to \mathbb R$ a function of class $C^2$ such that $\inf
f>-\infty$ and $\sup ||Hess f||<\infty$. Then every minimizing
sequence of $f$ admits a good shadow.
\end{thm}

A compelling reason for examining the above conjecture -- one that
goes beyond a mere comparison between the Ekeland variational
principle and the Yau minimum principle -- is that the good shadow
property is a useful tool to have if one wants to refine the
applications of the asymptotic minimum principle, as it allows for
information to be localized at infinity. Indeed, this is the
philosophy behind the proofs of Theorems \ref{ConvHull} and
\ref{ConvHullC2} below.

The remarks of the last paragraph are best understood in the
applications of the asymptotic minimum (maximum) principle to the
study of submanifolds. In these problems, geometric intuition can
often help to locate in space a particular (but,  \it a priori\rm,
not sufficiently well-behaved) minimizing sequence. One can then
try to find a good shadow of this sequence, for which the
computations of the relevant quantities can yield the desired
result.

In what follows we describe how the ideas outlined above can be
used, in the context of Theorem \ref{SS}, to study the problem of
characterizing the convex hulls of $C^1$ immersed submanifolds
that have  controlled Gauss maps.

Due to the low regularity, the usual tools of submanifold
geometry, centered as they are on the study of the second
fundamental form, cannot be applied to  $C^1$ submanifolds.
Nevertheless, these are naturally occurring objects, worthy of
study. For instance, it follows from a theorem of Nash-Kuiper
\cite{K} that every Riemannian manifold $M^m$ admits an isometric
$C^1$-embedding into an arbitrarily small neighborhood of
$R^{2m}$. Of course, for smoother immersions the expected
codimension is much higher \cite{G}. See also the end of this
Introduction for an interesting question regarding $C^1$ isometric
immersions, in the context of the present paper.

Any non-empty open convex subset $\cal O$ of $ \mathbb R^n$ is the
convex  hull of a $C^{\infty}$ complete submanifold, of any
codimension. To see this when $n\geq 3$, take a smooth curve
$\Gamma\subset \cal O$, of infinite length on both ends, whose
convex hull is $\cal O$.  Let $M$ be the union over all $p\in
\Gamma$ of smoothly varying $k$-dimensional spheres
$S^{(k)}_{r(p)}$, $1\leq k\leq n-2$, centered at $p$ and contained
in the normal space of $\Gamma$ at $p$.  Taking $r(p)$ to decay
fast enough   one can make sure that the resulting manifold $M$,
which is automatically complete, is contained in $ \cal O$.  Since
$p$ is in the convex hull of  $S^{(k)}_{r(p)}$ for any $p\in
\Gamma$, it follows that the convex hull of $M$ satisfies
$\text{Conv}\, (M)=\cal O$.

Using Theorem \ref{SS}, we show  that there are obstructions for a
given convex set to be the convex hull of a complete $C^1$
submanifold of a fixed codimension, provided the  Gauss map, which
is of course continuous, is assumed to be  {\it uniformly}
continuous. \rm

Recall that  $h:M^m\to \mathbb R^n$ is said to be \it substantial
\rm if $h(M)$ is not contained in a proper affine subspace of
$\mathbb R^n$.

\begin{thm}\label{ConvHull} Let $M$ be a complete $m$-dimensional Riemannian manifold, $n>m$,
and $h:M \to\mathbb R^n$ a substantial $C^1$ isometric immersion
for which the Grassmanian-valued Gauss map $\mathcal G:M\to G(n-m,
n)$, given by $\mathcal G(p)= [h_{\ast}TM_p]^{\perp}$, is
uniformly continuous. Then either $\text{Conv}\, [h(M)]=\mathbb
R^n$, or each point in the boundary of $\text{Conv}\, [h(M)]$
admits at most $n-m$ supporting hyperplanes in general position.
\end{thm}

\begin{cor}\label{Comp} If $h: M^m \to \mathbb R^n$ is a substantial $C^1$
immersion of a compact manifold, then each point in the boundary
of $\text{Conv}\, [h(M^m)]$ admits at most $n-m$ supporting
hyperplanes in general position.
\end{cor}

Observing that the limit of supporting hyperplanes is itself a
supporting hyperplane, we have:

\begin{cor}\label{Comp1}
If $M^{n-1}$ is compact and $h:M^{n-1}\to\mathbb R^n$ is a $C^1$
immersion, then $h$ is substantial and each point $p$ in the
boundary of $\text{Conv}\,[h(M^{n-1})]$ admits a unique support
hyperplane $H_p$. Moreover, the map $p\mapsto H_p$ is continuous.
\end{cor}

Theorem \ref{weakyau} also has applications to the geometry of
submanifolds:

\begin{thm}\label{ConvHullC2} Let $M^m$ be a complete manifold, $n>m$,
and $h: M^m\to\mathbb R^n$ a substantial $C^2$ isometric
immersion, with bounded second fundamental form and uniformly
continuous mean curvature vector field $\overrightarrow H$.
Suppose $\text{Conv}\,[h(M^m)]\neq\mathbb R^n$ and let
$H_1,\ldots,H_s$ be supporting hyperplanes in general position
that pass through a point $p_o$ of the boundary of
$\text{Conv}\,[h(M^m)]$. Let $e_i\in [H_i]^{\perp}$ be the unit
vector such that $\langle h(x)-p_o, e_i\rangle\geq 0$ for every
$x\in M^m,\;i=1,...,s$. Then $s\leq n-m$ and there exists a
sequence $(p_k)$ in $M^m$ such that $\text{d}(h(p_k),
H_1\cap\dots\cap H_s )\to 0$ and
$\liminf_{k\to\infty}\left\langle\overrightarrow H(p_k),e_i
\right\rangle\geq 0,\;i=1,\ldots,s$.
\end{thm}

The following result can be viewed as a generalization of the fact
that the mean curvature of a compact convex hypersurface of
$\mathbb R^n$ is nonnegative, after an appropriate choice of the
orientation.

\begin{thm}\label{ConvHullC2a} Let $h:M^m\to\mathbb R^n$ be a substantial $C^2$ isometric immersion
of a compact Riemannian manifold $M^m$, and $H_1,\ldots,H_s$
supporting hyperplanes in general position that pass through a
point $p_o$ of the boundary of $\text{Conv}\,[h(M^m)]$. Let
$e_1,\ldots,e_s$ be as in the statement of the Theorem
\ref{ConvHullC2}. Then $s\leq n-m$, $h(M)\cap H_1\cap\dots\cap
H_s\neq\emptyset$, and for every point $q\in M^n$ such that
$h(q)\in H_1\cap\dots\cap H_s$ one has
$\left\langle\overrightarrow H(q),e_i\right\rangle\geq
0,\;i=1,\ldots,s$.
\end{thm}

\begin{figure}[hbt]
\centerline{\resizebox{7cm}{!}{\includegraphics{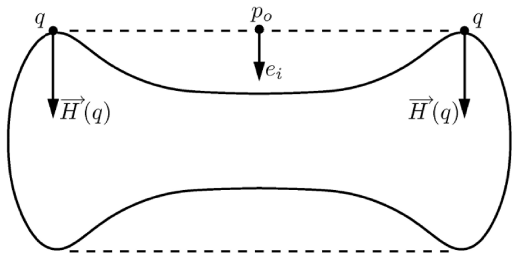}}}
\caption{Theorem 1.8}
\end{figure}

\vskip5pt

\noindent \bf Examples. \rm  It is easy to illustrate Theorem
\ref{ConvHull}, already in low dimensions:

\vskip10pt

\noindent i)  Let $l$ be a line in $\mathbb R^3$, and $P_1$, $P_2$
planes such that $P_1\cap P_2=l$. Let $\cal O$ be a component of
$\Bbb R^3-(P_1\cup P_2)$. One can construct a complete
$C^{\infty}$ curve $\Gamma \subset \cal O$ such that $\text{Conv}
\;(\Gamma)=\cal O$ and $\Gamma$ has bounded curvature. The last
condition ensures that the Gauss map $\mathcal G:\Gamma \to
G(2,3)$ is uniformly continuous. Along $l$, the maximum number of
supporting hyperplanes to $\partial \cal O$ that are in general
position is two, which is also the codimension of $\Gamma$.
This gives the equality case in Theorem \ref{ConvHull}.

We give an informal description of how $\Gamma$ can be
constructed. Start with oriented line segments $l_n$ parallel to
$l$, $n\geq 1$, contained in $\cal O$, getting longer as $n\to
\infty$, and accumulating onto the entire oriented line $l$.  One
obtains $\Gamma$ by  connecting  for all $n\geq 1$   the last
point of $l_n$, in a smooth way,  to the first point of $l_{n+1}$,
by means of a curve $\gamma_n$  of curvature less than one. The
curve $\gamma_n$ is supposed to be very long, going deep inside
$\cal O$ and turning slowly, so that the curvature can be kept
smaller than one. Once $\gamma_n$ is far  from $l$,  one can also make $\gamma_n$
twist around, with controlled curvature,  so as to make its convex
hull bigger.  It is now clear that a sequence of curves $\gamma_n$
can be created so that  $\Gamma$ has  curvature less than one and
$\text{Conv} \;(\Gamma)=\cal O$.

Observe that such a construction is impossible if, instead of a
curve, one takes $\Gamma$ to be a complete surface. Indeed, as the
surface gets closer and closer to $l$, in order for $\Gamma$ to
remain in $\cal O$ it has to fold abruptly, thus violating the
condition that the Gauss map is uniformly continuous.

\vskip10pt

\noindent ii) Let $M$ be an open hemisphere in $S^2\subset \mathbb
R^3$. Its  convex hull  is, of course,  the solid hemisphere. At
points along the great circle, $\text{Conv}\; (M)$ has two
supporting hyperplanes in general position,  whereas the
codimension of $M$ is one.  This shows that Theorem \ref{ConvHull}
fails if the submanifold is not complete.\qed

\vskip10pt

The uniform continuity condition on the Gauss map allows for the
Ricci curvature of the submanifold to be unbounded from below. In
fact, it is easy to construct smooth complete graphs in $\mathbb
R^3$ with these properties. This shows that  Theorem \ref{OY}
cannot be applied to prove Theorem \ref{ConvHull}, even if the
submanifold in question is of class $C^{\infty}$.

To put these remarks in perspective note that, by the Gauss
equation, the natural way to force the intrinsic curvatures of a
submanifold to be bounded is simply to require that the second
fundamental form has bounded length. Although this is not obvious,
at least in the case of hypersurfaces the latter condition  means
that the Gauss map is globally Lipschitzian, which is stronger
than merely requiring the Gauss map to be uniformly continuous.

We stress that, in Theorem \ref{ConvHull}, even if the submanifold
is $C^{\infty}$ and has bounded second fundamental form, the Yau
minimum principle cannot be applied. Indeed, as it will be clear
from the proof, one needs to find good shadows that are provided
in the $C^1$ context by Theorem \ref{SS},  for \it arbitrary \rm
minimizing sequences.  Theorem \ref{OY}, on the other hand,
guarantees the existence of a \it single \rm minimizing sequence
with \lq\lq good\rq\rq properties.

A natural question is whether the condition in Theorem
\ref{ConvHull}, stating that the boundary points of the convex
hull admits at most $n-m$ tangent hyperplanes in general position,
is also sufficient for the construction of examples. Here, again,
one ought to keep in mind that the Gauss map is only required to
be uniformly continuous, instead of the stronger condition of
being Lipschitz. We are indebted to J. Fu for pointing out that
the work of Alberti \cite{A} may be relevant to this question.

Another interesting problem concerning immersions with uniformly
continuous Gauss maps, albeit unrelated to the actual results in
the present paper, has to do with the isometric immersions
$\varphi:M^2\to R^3$ of Hadamard surfaces with curvature bounded
away from zero. By Efimov's theorem \cite{Ef}, no such $\varphi$
exists that is of class $C^2$. On the other hand, by the
Nash-Kuiper's theorem \cite{K} there is a $\varphi$ of class
$C^1$; in particular, its Gauss map is continuous. If $M^2$ is as
above, does there exist a $C^1$ isometric immersion
$\varphi:M^2\to \mathbb R^3$ whose Gauss map is uniformly
continuous? By \cite{K}, the answer is yes if $M^2$ is the
universal cover of a compact surface with negative curvature.

Finally, we refer to \cite{X} for other geometric properties that
can be recovered from Grassmanian-valued Gauss maps.

\newpage

\section{A strong {$\bf C^1$} version of the Yau minimum principle}

\noindent For the first proof of Theorem \ref{SS} we will need the
following fundamental result (\cite{E}, \cite{Fi}, \cite{St}):

\vskip10pt

\noindent{\bf The Ekeland Variational Principle.} {\it Let $(X,d)$
be a complete metric space, and $f:X\to \Bbb R$ a function which
is lower semi-continuous and bounded from below. Then for any
$\varepsilon, \delta>0$, and $x\in X$ with
$f(x)\leq\inf_Xf+\varepsilon$, there is $y\in X$ satisfying}

\noindent{i)} $d(x,y)\leq\delta$

\noindent{ii)} $f(y)\leq f(x)$

\noindent{iii)}
$f(y)<f(z)+\frac{\varepsilon}{\delta}d(y,z)$,\;\;\;for all $z\in
X$ with $z\neq y$.

\vskip15pt

\noindent{\bf First proof of Theorem \ref{SS}.} For each
$n\in\mathbb N$, let $\varepsilon_n=f(x_n)-\inf_Mf$,
$\delta_n=\sqrt{\varepsilon_n}$. In the sequel we will prove the
existence of a sequence $(y_n)$ in $M$ satisfying
\begin{eqnarray}\label{SSa}
f(y_n)\leq f(x_n),\;\;\; d(x_n,y_n)\leq \delta_n
\end{eqnarray}
and
\begin{eqnarray}\label{SSb}
||\nabla f||(y_n)\leq \delta_n.
\end{eqnarray}
Since $\delta_n\to 0$ as $n\to\infty$, $(y_n)$ will have the
desired properties. If $f(x_n)=\inf_Mf$, take $y_n=x_n$.
Otherwise, we have $\varepsilon_n>0$, $\delta_n>0$, and applying
the Ekeland Variational Principle with
$\varepsilon=\varepsilon_n$, $\delta=\delta_n$ and $x=x_n$, we
obtain $y_n\in M$ satisfying (\ref{SSa}) and
\begin{eqnarray}\label{SS3}
f(y_n)<f(z)+\frac{\varepsilon_n}{\delta_n}d(y_n,z)=f(z)+\delta_n\,d(y_n,z),
\end{eqnarray}
for all $z\in M$ with $z\neq y_n$. To show (\ref{SSb}), take an
arbitrary unit vector $v\in T_{y_n}M$ and let $\gamma : (-c,c)\to
M$ be the unit speed geodesic in $M$ so that $\gamma(0)=y_n$ and
$\gamma'(0)=v$. Reducing $c$ if necessary, we can suppose that the
image of $\gamma$ is contained in a normal neighborhood of $y_n$
in M. From (\ref{SS3}) we have, with $z=\gamma(t)$,
\begin{eqnarray}\label{SS4}
f\big( \gamma(t)\big)-f(y_n)>-\delta_n\,d\big(
\gamma(t),y_n\big)=-\delta_n\,|t|,\;\;\;0<|t|<c,
\end{eqnarray}
which implies
\begin{eqnarray}\label{SS5}
\frac{f\big( \gamma(t)\big)-f(y_n)}{t}<\delta_n,\;\;\;-c<t<0.
\end{eqnarray}
Since $f$ is of class $C^1$, it follows that
\begin{eqnarray}\label{SS6}
\langle\nabla
f(y_n),v\rangle=\frac{d}{dt}\Big|_{t=0}f\circ\gamma(t)=\lim_{t\to
0^-}\frac{f\big( \gamma(t)\big)-f(y_n)}{t}\leq\delta_n,
\end{eqnarray}
for all $v\in T_{y_n}M$ with $||v||=1$. Therefore,
\begin{eqnarray}\label{SS7}
\big| \langle \nabla f(y_n),v\rangle\big|\leq\delta_n,
\end{eqnarray}
for all unit vector $v\in T_{y_n}M$, so that
\begin{eqnarray}\label{SS7}
||\nabla f(y_n)||\leq\delta_n, \;\;\;n\in\mathbb N.
\end{eqnarray}
The sequence $(y_n)$ satisfies (\ref{SSa}) and (\ref{SSb}) and
thus the conditions of the theorem. \qed

\vskip15pt

\noindent{\bf Second proof of Theorem  \ref{SS}.} For each
$n\in\mathbb N$, let
\begin{eqnarray}\label{SP1}
r_n=\sqrt{f(x_n)-\inf_Mf}.
\end{eqnarray}
As in the first proof of Theorem \ref{SS}, we will construct a
sequence $(y_n)$ in $M$ satisfying
\begin{eqnarray}\label{SPnew}
f(y_n)\leq f(x_n),\;\;\; d(x_n,y_n)\leq r_n,\;\;\; ||\nabla
f||(y_n)\leq r_n,
\end{eqnarray}
for all $n\in\mathbb N$. The conclusion that $(y_n)$ meets the
conditions of the theorem will follow from the fact that $r_n\to
0$ when $n\to\infty$. If $\nabla f(x_n)=0$, take $y_n=x_n$. In
case $\nabla f(x_n)\neq 0$, denote by $\gamma_n:[0, \tau_n)\to M$
the maximal integral curve of the vector field $X=-||\nabla
f||^{-2}\, \nabla f$ such that $\gamma_n(0)=x_n$. From $
(f\circ\gamma_n)'=-1$,  we obtain
\begin{eqnarray}\label{SP2}
f(\gamma_n(t))-f(x_n)=-t
\end{eqnarray}
for all $t$ in $[0, \tau_n)$. Recall that $q\in M$ belongs to the
$w$-limit set of $\gamma_n$ if there exists a sequence $(t_k)$
contained in $[0, \tau_n)$ such that $\lim_{k\to \infty}
t_k=\tau_n$ and $\lim_{k\to \infty} \gamma_n (t_k)=q$. The orbits
of X are also orbits of $-\nabla f$. It follows from (\cite{PM},
p. 13) that the $w$-limit set of $\gamma_n$ consists entirely of
critical points of $f$. Denote by $\overline B(x_n, r_n)$ the
closed ball with center at $x_n$ and radius $r_n$. We have two
cases to consider:

\vskip5pt

\noindent{i)} The $\omega$-limit set of $\gamma_n$  contains a
point $q_n$ in $\overline B(x_n,r_n)$. As observed above, $q_n$ is
necessarily a critical point for $f$ and we set $y_n=q_n$.

\vskip5pt

\noindent{ii)} The $\omega$-limit set of $\gamma_n$ does not
intersect $\overline B(x_n,r_n)$. Let $t_n$ be the first time such
that the positive trajectory of $X$ through $x_n$ hits $\partial
B(x_n, r_n)$. We have
\begin{eqnarray}\label{SP3}
r_n \leq \hbox{length} \;\gamma_n|_{[0, t_n]}= \int_0^{t_n}
||\gamma_n'(t)||dt\leq  t_n \; \max _{\gamma_n|_{[0,t_n]}}
(||\nabla f||^{-1}),
\end{eqnarray}
so that, by (\ref{SP2}),
\begin{eqnarray}\label{SP4}
r_n \min_{\gamma_n|_{[0,t_n]}} ||\nabla f||\leq t_n=
f(x_n)-f(\gamma_n(t_n))< f(x_n)-\inf_Mf.
\end{eqnarray}
It follows from (\ref{SP1}) and (\ref{SP4}) that
\begin{eqnarray}\label{SP5}
\min_{\gamma_n|_{[0,t_n]}} ||\nabla f||<r_n.
\end{eqnarray}
Choosing $\theta_n\in [0,t_n]$ so that $||\nabla
f||(\gamma_n(\theta_n))=\min_{\gamma_n|_{[0,t_n]}} ||\nabla f||$
and taking $y_n=\gamma_n(\theta_n)$, we have $f(y_n)\leq f(x_n)$
and, since $\gamma_n([0,t_n])\subset\overline B(x_n,r_n)$,
$$
d(x_n,y_n)\leq r_n=\sqrt{f(x_n)-\inf_Mf}.
$$
Moreover, by (\ref{SP5}), $||\nabla f||(y_n)<r_n$. Hence the
sequence $(y_n)$ satisfies (\ref{SPnew}) and thus the conditions
of the theorem.\qed

\section {Flows and asymptotic minimum principles}

In this section we study the conjecture from the Introduction,
using gradient flows as in the second proof of Theorem (\ref{SS}).
Here, however, the situation is much more complex and only partial
results are available. Instead of estimating only lengths, one
needs to estimate both lengths and volumes. Further refinements of
the basic strategy may yet yield a proof of the full conjecture.

\vskip10pt

\noindent{\bf Proof of Theorem \ref{weakyau}.} Let $\phi_t$ be the
local flow of $X=-\nabla f$ on $M$, so that
\begin{eqnarray} \frac {d}{dt}\phi_t(p)=-\nabla f(\phi_t(p)), \;\;
\phi_0(p)=p, \label{ode}
\end{eqnarray}
where  $t\in [0,\tau(p))=$ the maximal interval of existence of
the forward solution.

For each $n\in\mathbb N$, set
\begin{eqnarray}\label{weakyau1}
\delta_n=\sqrt{f(p_n)-\inf_Mf},\;\;\;r_n=\sqrt{\delta_n}.
\end{eqnarray}
Since $\delta_n\to 0$ as $n\to\infty$, we may suppose, without
loss of generality, that $\delta_n<r_n$ for all $n\in\mathbb N$.
We will construct a sequence $(q_n)$ in $M$ satisfying
\begin{eqnarray}\label{weakyau2}
\text{d}(p_n,q_n)\leq r_n,\;\;\; \liminf_{n\to\infty}\Delta
f(q_n)\geq 0.
\end{eqnarray}
If $\delta_n=0$, we have $\Delta f(p_n)\geq 0$ and we choose
$q_n=p_n$. Fix $n\in\mathbb N$ so that $\delta_n>0$. We are going
to distinguish between two cases:

\vskip5pt

\noindent{a)} Every positive orbit originating in $\overline
B(p_n,r_n^2)$ remains in the open ball $B(p_n,r_n)$.

\vskip5pt

\noindent{b)} There is at least one trajectory that joins the
boundaries of $B(p_n,r_n^2)$ and $B(p_n,r_n)$ in finite time.

In the first alternative, $\tau(p)=\infty$ for every
$p\in\overline B(p_n,r_n^2)$. Let $\mu$ denote the Riemannian
measure of $M$. By Liouville's formula for the change of volume
under a flow \cite{M}, one has, for all $t>0$, since
$\Delta=\text{div}\nabla$,
\begin{eqnarray}\label{weakyau3}
\mu\big(B(p_n,r_n)\big)\geq \mu
\big(\phi_t(B(p_n,r_n^2))\big)=\int_{B(p_n,r_n^2)} \exp \Big(
\int_0^t -\Delta f (\phi_s(p))ds\Big)d\mu(p).
\end{eqnarray}
If there exists $\varepsilon >0$ such that $\Delta f(q)\leq
-\varepsilon$ for all $q\in B(p_n,r_n)$, a contradiction can be
easily established by letting $t\to \infty$ in the above formula.
Hence one can choose $q_n\in B(p_n,r_n)$ such that $\Delta f
(q_n)\geq 0$.

We now work under the conditions of alternative b). Consider the
quantity $\tau_n$ which gives the shortest time to travel from
$\partial B(p_n,r_n^2)$ to $\partial B(p_n, r_n)$, along a
trajectory of $X$. Formally,
\begin{eqnarray}\label{weakyau4}
\tau_n= \inf \;\big\{t: t\in (0, \tau(p)),\;p\in\partial
B(p_n,r_n^2)\; \text{and} \;\phi_t(p)\in \partial B(p_n, r_n)
\big\}.
\end{eqnarray}
In particular,
\begin{eqnarray}\label{weakyau5}
\phi_{\tau_n}\big(\,\overline B(p_n,r_n^2)\big)\subset\overline
B(p_n, r_n).
\end{eqnarray}

We want to estimate the first exit time $\tau_n$. Let
$x_n\in\partial B(p_n,r_n^2)$ and $t_n\in (0,\tau(x_n))$ be such
that $\phi_{t_n}(x_n)\in\partial B(p_n,r_n)$ and $t_n<2\tau_n$.

The last integral of
\begin{eqnarray}\label{weakyau6}
f(x_n)- f(\phi_{t_n}(x_n))= - \int_0^{t_n} (f\circ \phi_s)'ds&=&
\int_0^{t_n}||\nabla f||^2 (\phi_s(x_n))ds\nonumber\\ &\geq&
\frac{1}{t_n}\left [ \int_0^{t_n}||\nabla f|| (\phi_s(x_n))ds
\right]^2
\end{eqnarray}
gives the length of the portion of the orbit of $X=-\nabla f$ over
the time interval $[0,t_n]$, through $x_n$. Since the last point
of this orbit segment lies in $\partial B(p_n,r_n)$, its length is
at least $r_n(1-r_n)$. Collecting this information, and observing
(\ref{weakyau6}), we obtain

\begin{eqnarray}\label{weakyau7}
\frac{\delta_n(1-r_n)^2}{2\tau_n}\leq\frac{r_n^2(1-r_n)^2}{t_n}
\leq f(x_n)- f(\phi_{t_n}(x_n))\leq f(x_n) -\inf_Mf.
\end{eqnarray}

We will now estimate the last term in (\ref{weakyau7}). Define
$h:M\to\mathbb R$ by $h(x)=||\nabla f(x)||^2 +\varepsilon$, where
$\varepsilon$ is a positive real number. Given $p, q\in M$,
consider an unit speed minimizing geodesic $\gamma:[0,a]\to M$
joining $p$ to $q$. If $K$ is an upper bound for the norm of the
Hessian operator of $f$, we have
\begin{eqnarray}\label{lemweakyau2}
\Big|\frac{d}{dt} h(\gamma(t))\Big|=2|\langle
\nabla_{\gamma'}\nabla f, \nabla f\rangle|\leq 2K||\nabla f||\leq
2K\sqrt{h(\gamma(t))},
\end{eqnarray}
and so
\begin{eqnarray}\label{lemweakyau3}
\big|\sqrt{h(\gamma(t))}-\sqrt{h(\gamma(0))}\,\big|\leq Kt,\;\;\;
t>0.
\end{eqnarray}
Setting $t=a$ and letting $\varepsilon\to 0$,
\begin{eqnarray}\label{lemweakyau4}
\big|\,||\nabla f(p)||-||\nabla f(q)||\,\big|\leq
Ka=Kd(p,q),\;\;\;p, q\in M.
\end{eqnarray}
On the other hand, from the second proof of Theorem \ref{SS},
there exists $y_n\in\overline B(p_n,\delta_n)$ so that $||\nabla
f||(y_n)\leq\delta_n$. Using this fact and (\ref{lemweakyau4}), we
obtain
\begin{eqnarray}\label{weakyau8}
||\nabla f(z)||\leq ||\nabla
f(y_n)||+K\text{d}(y_n,z)\leq\delta_n(1+2K),\;\;\;z\in \overline
B(p_n,\delta_n).
\end{eqnarray}

Considering an unit speed minimizing geodesic segment
$\gamma:[0,\delta_n]\to M$ joining $p_n$ to $x_n$, it follows from
(\ref{weakyau8}) that
\begin{eqnarray}\label{weakyau9}
|f(x_n)-f(p_n)|&\leq&\int_0^{\delta_n}|(f\circ\gamma)'(t)|dt=\int_0^{\delta_n}|\langle
\nabla
f(\gamma(t)),\gamma'(t)\rangle|dt\nonumber\\&\leq&\int_0^{\delta_n}||\nabla
f(\gamma(t))||dt\leq\delta_n^2(1+2K).
\end{eqnarray}
From (\ref{weakyau7}) and (\ref{weakyau9}), we obtain
\begin{eqnarray}\label{weakyau10}
\frac{\delta_n(1-r_n)^2}{2\tau_n}&\leq&
f(x_n)-f(p_n)+f(p_n)-\inf_Mf\nonumber\\&\leq&\delta_n^2(1+2K)+\delta_n^2=2\delta_n^2(1+K),
\end{eqnarray}
and so
\begin{eqnarray}\label{weakyau11}
\frac{(1-r_n)^2}{\tau_n}\leq 4\delta_n(1+K).
\end{eqnarray}

\noindent{Since} $\delta_n\to 0$ and $r_n\to 0$ as $n\to\infty$,
one has, in particular, $\lim_{n\to\infty}\tau_n=+\infty$.

\vskip5pt

From Liouville's formula,

\begin{eqnarray}\label{weakyau12}
\mu\big(\phi_{\tau_n}(B(p_n,r_n^2))\big)=\int_{B(p_n,r_n^2)}\exp\left(\int_0^{\tau_n}
-\Delta f (\phi_s(p))ds\right)d\mu(p)\label {fctn1}.
\end{eqnarray}

Jensen's inequality applied to the probability measure $\nu/
\nu(\Omega)$, where $\nu$ is a finite measure on $\Omega$, gives
\begin{eqnarray}\label{weakyau13}
\psi\left (\frac{1}{\nu(\Omega)}\int_{\Omega} gd\nu\right)\leq
\frac{1}{\nu(\Omega)}\int_{\Omega}(\psi\circ g )  d\nu,
\end{eqnarray}
whenever $\psi$ is convex and $g$ is integrable.

\vskip5pt

Applying (\ref{weakyau13})  to (\ref{weakyau12}),

\begin{eqnarray}\label{weakyau14}
\frac {\mu(\phi_{\tau_n}(B(p_n,r_n^2)))}{\mu(B(p_n,r_n^2))}\geq
\exp \left [ \frac{1}{ \mu(B(p_n,r_n^2))}
\int_{B(p_n,r_n^2)}\left(\int_0^{\tau_n} -\Delta
f(\phi_s(p))ds\right)d\mu(p)\right].
\end{eqnarray}

\vskip10pt

If $M^c$ stands for the space form of curvature $c<0$, where $c$
is sufficiently negative as compared to the lower bound of the
Ricci curvature of $M$, Gromov's theorem on monotonicity of volume
ratios (\cite{CH}, p. 125) implies that for all $x\in M$ the
quotient
$$\frac {v^c(r)}{v(x,r)},$$
between the volumes of balls of radius $r$ in $M^c$ and $M$, is a
nondecreasing function. In particular,
\begin{eqnarray}\label{weakyau15}
\frac{\mu\big(B(p_n,r_n)\big)}{\mu\big(B(p_n,r_n^2)\big)}\leq\frac{v^c(r_n)}{v^c(r_n^2)}.
\end{eqnarray}

\noindent{Recalling} that $\phi_{\tau_n}(\overline
B(p_n,r_n^2))\subset\overline B(p_n,r_n)$, it follows from
(\ref{weakyau14}) and (\ref{weakyau15}) that
\begin{eqnarray}\label{weakyau16}
\frac{1}{\mu(B(p_n,r_n^2))}\int_{B(p_n,r_n^2)}
\Big(\int_0^{\tau_n}-\Delta f(\phi_s(p))ds\Big) d\mu(p)\leq
\log\Big(\frac{v^c(r_n)}{v^c(r_n^2)}\Big).
\end{eqnarray}

Next, we take $\Omega_n= B(p_n,r_n^2)\times [0,\tau_n]$, endowed
with the probability measure $\nu$ given by the normalization of
the product measure on $B(p_n,r_n^2)\times [0,\tau_n]$. Dividing
(\ref{weakyau16}) by $\tau_n$, and using (\ref{weakyau11}), we
arrive at
\begin{eqnarray}\label{weakyau18}
\int_{\Omega_n} - \Delta f(\phi_s(q))d\nu(q,s)
\leq\frac{1}{\tau_n}\log\Big(\frac{v^c(r_n)}{v^c(r_n^2)}\Big)&\leq&\frac{4(1+K)}{(1-r_n)^2}\delta_n
\log\Big(\frac{v^c(r_n)}{v^c(r_n^2)}\Big)\nonumber\\&\leq& 16
(1+K)\delta_n \log\Big(\frac{v^c(r_n)}{v^c(r_n^2)}\Big),
\end{eqnarray}
for every sufficiently large $n$.

Since  $\Omega_n$ has mass one, it follows from (\ref{weakyau18})
that, for every sufficiently large $n\in \mathbb N$, there are
$s_n\in [0,\tau_n]$ and $q_n'\in B(p_n,r_n^2)$ such that, with
$q_n=\phi_{s_n}(q_n')$, one has
\begin{eqnarray}\label{weakyau19}
d(p_n, q_n)\leq r_n, \;\;\;\Delta f(q_n)\geq
16(1+K)\delta_n\log\Big(\frac{v^c(r_n^2)}{v^c(r_n)}\Big).
\end{eqnarray}

The volume element of $M^c$ in spherical coordinates is uniformly
bounded from above and below by fixed multiples of $r^{m-1}$ if
$r<1$. In particular, there exists a constant $A>0$ such that
\begin{eqnarray}\label{weakyau19a}
\delta_n\log\Big(\frac{v^c(r_n^2)}{v^c(r_n)}\Big)\geq\delta_n\log
\Big(A \frac{r_n^{2m}}{r_n^m}\Big)=\delta_n\log(A\,r_n^m),
\end{eqnarray}
which implies, since $\delta_n=r_n^2\to 0$ as $n\to\infty$,
\begin{eqnarray}\label{weakyau20}
\liminf_{n\to\infty}\Delta f(q_n)\geq
16(1+K)\lim_{n\to\infty}\delta_n\text{log}(\delta_n^{m/2})
=8(1+K)m\lim_{n\to\infty}\delta_n\text{log}\,\delta_n=0.
\end{eqnarray}
Thus the sequence $(q_n)$ satisfies (\ref{weakyau2}) as desired.
To complete the proof of the theorem it remains to show that
$||\nabla f||(q_n)\to 0$ as $n\to\infty$. By Theorem \ref{SS},
there exists a minimizing sequence $(q_n')$ with $\text{d}(q_n,
q_n')\to 0$, $||\nabla f||(q_n')\to 0$. Applying
(\ref{lemweakyau4}) to $q_n$ and $q_n'$, one sees that $||\nabla
f||(q_n)$ also tends to zero. \qed

\vskip10pt

Adjusting the proof of Theorem 1.3 one has the following:

\vskip10pt

\begin{thm}\label{oscil} Let $M$ be a complete manifold with Ricci curvature bounded from below,
and $f:M\to \mathbb R$ a function of class $C^2$ satisfying
$\inf_Mf>-\infty$. Let $(p_n)$ be a sequence in $M$ that is
strongly minimizing for $f$, in the sense that there exists
$\delta>0$ such that the oscillation of $f$ on $B(p_n, \delta)$
tends to zero, i.e.,
\begin{eqnarray}\label{oscil1}
\lim_{n\to \infty} \left [ \max_{ B(p_n, \delta)}f - \min_ {
B(p_n, \delta)}f \right ]=0.
\end{eqnarray}
Then there exists a minimizing sequence $(q_n)$ in $M$ for $f$
such that
\begin{eqnarray}\label{oscil2}
\lim_{n\to\infty}d(p_n,q_n)=0,\;\;\;\liminf_{n\to \infty} \Delta
f(q_n)\geq 0.
\end{eqnarray}
\end{thm}

\noindent{\bf Proof.} For each $n\in\mathbb N$, set
\begin{eqnarray}\label{oscil3}
r_n=\left [ \max_{ B(p_n, \delta)}f - \min_ { B(p_n, \delta)}f
\right ]^{\frac{1}{4}}.
\end{eqnarray}
We will construct a sequence $q_n$ such that
\begin{eqnarray}\label{oscil4}
\text{d}(p_n,q_n)\leq r_n,\;\;\;\liminf\Delta f(q_n)\geq 0.
\end{eqnarray}
From (\ref{oscil1}), one sees that such a sequence $(q_n)$ will
satisfy (\ref{oscil2}). Since $r_n\to 0$ as $n\to\infty$, we may
suppose, without loss of generality, that $r_n\in[0,\delta)$ for
all $n\in\mathbb N$, where $\delta$ is as in the statement of the
theorem. If $r_n=0$, $f$ is constant in $B(p_n,\delta)$, and we
take $q_n=p_n$. Fix a real number $\kappa\in (0,1)$ and denote by
$\phi_t$ the local flow of $X=-\nabla f$ on $M$. For each
$n\in\mathbb N$ for which $r_n>0$, we have two possibilities:
either every positive orbit originating in $\overline B(p_n,\kappa
r_n)$ remains in the open ball $B(p_n,r_n)$ or there is at least
one trajectory that joins the boundaries of $B(p_n,\kappa r_n)$
and $B(p_n,r_n)$ in finite time. In the first case we prove, in
the same way as in the proof of Theorem (\ref{weakyau}), the
existence of $q_n\in\overline B(p_n,r_n)$ such that $\Delta
f(q_n)\geq 0$, and we take $q_n=p_n$. In the second case,
reasoning as in the proof of Theorem \ref{weakyau}, with $r_n^2$
replaced by $\kappa r_n$, we obtain
\begin{eqnarray}\label{oscil5}
\frac{(1-\kappa)^2r_n^2}{2\tau_n}\leq\max_{B(p_n,\delta)}f-\min_{B(p_n,\delta)}f,
\end{eqnarray}
which implies, in view of (\ref{oscil3}),
\begin{eqnarray}
\frac{(1-\kappa)^2}{2\tau_n}\leq r_n^2.
\end{eqnarray}
In particular, $\tau_n\to\infty$ as $n\to\infty$. Continuing as in
the proof of Theorem \ref{weakyau}, we arrive at
\begin{eqnarray}\label{oscil6}
\frac {\mu(B(p_n,r_n))}{\mu(B(p_n,\kappa r_n))}\geq \exp \left [
\frac{1}{ \mu(B(p_n,\kappa r_n))} \int_{B(p_n,\kappa
r_n)}\left(\int_0^{\tau_n} -\Delta
f(\phi_s(p))ds\right)d\mu(p)\right].
\end{eqnarray}
By Gromov's theorem on monotonicity of volume ratios (\cite{CH},
p. 125),
\begin{eqnarray}\label{oscil7}
\frac{\mu\big(B(p_n,r_n)\big)}{\mu\big(B(p_n,\kappa
r_n)\big)}\leq\frac{v^c(r_n)}{v^c(\kappa r_n)},
\end{eqnarray}
where $c$ is sufficiently negative as compared to the lower bound
of the Ricci curvature of $M$, and $v^c(r)$ is the volume of a
closed ball in the space form $M^c$ of curvature $c<0$. Since
$r_n<\delta$ for all $n$, it follows from the homogeneity of $M^c$
that there exists $C>0$ such that
\begin{eqnarray}\label{oscil8}
\frac{v^c(r_n)}{v^c(\kappa r_n)}\leq C.
\end{eqnarray}
It follows from (\ref{oscil6}), (\ref{oscil7}) and (\ref{oscil8})
that
\begin{eqnarray}\label{oscil9}
\frac{1}{\mu(B(p_n,\kappa r_n))}\int_{B(p_n,\kappa r_n)}
\Big(\int_0^{\tau_n}-\Delta f(\phi_s(p))ds\Big) d\mu(p)\leq \log
C.
\end{eqnarray}
Using (\ref{oscil9}) and arguing as in Theorem \ref{weakyau}, we
conclude that there exists $q_n\in\overline B(p_n,r_n)$ so that
\begin{eqnarray}\label{oscil10}
\Delta f(q_n)\geq\frac{\text{log}\,(C^{-1})}{\tau_n}.
\end{eqnarray}
It is now clear that the sequence $(q_n)$ satisfies
(\ref{oscil4}), as desired.\qed

\vskip10pt

\noindent \bf Remarks. (i)\rm \,The conclusion in Theorem
\ref{oscil} fails if the Ricci curvature is unbounded from below.
To see this, let $X:M^2\to B(0,1)\subset R^3$ be a complete proper
minimal immersion. Examples of such surfaces were constructed by
Mart\'in-Morales \cite{MM}. It follows from \cite{X1} that the
Gaussian curvature of X is necessarily unbounded from below. Let
$f(p)=-|X(p)|^2$. Since X is proper in the unit ball, $\lim_{p\to
\infty} f(p)=-1$ uniformly on $p$, and condition (\ref{oscil1})
holds. On the other hand, minimality implies that $\Delta f=-4$,
so that the condition $\liminf_{n\to\infty} \Delta f(q_n)\geq 0$
can never be realized.

\vskip5pt

\noindent{\bf (ii)} Taking $\kappa=1/2$ one sees from
(\ref{oscil7}) that the hypothesis that the Ricci curvature is
bounded from below can be weakened to the condition that $M$
satisfies a {\it local volume doubling condition}: there exist $a,
b>0$ such that for any $p\in M$ and $0<r<a$, one has $\text{Vol}
\;B(p, r)\leq b \text{Vol}\;B(p,\frac{r}{2})$.

\vskip10pt

\section{Proofs of the geometric theorems}

\vskip10pt

\noindent{\bf Proof of Theorem \ref{ConvHull}.} Suppose
$\text{Conv}\,[h(M)]\neq\mathbb R^n$ and let $H_1, \cdots, H_s $
be supporting hyperplanes in general position through a point
$p_o$ in the boundary of $\text{Conv}\,[h(M)]$. We want to show
that $s\leq n-m$. To this end, for $i=1,...,s$, denote by $e_i$
the unit vector that is normal to $H_i$ and points inside
$\text{Conv}\,[h(M)]$, and let $f_i: \mathbb R^n\to\mathbb R$ be
the height function with respect to $H_i$, i.e.,
\begin{eqnarray}\label{heightfunction}
f_i(y)=\langle y-p_o,e_i\rangle.
\end{eqnarray}
The fact that $e_i$ points inside $\text{Conv}\,[h(M)]$,
$i=1,\ldots,s$, means that
\begin{eqnarray}\label{heightfunctionA}
h(M)\subset \{y\in\mathbb R^n:f_i(y)\geq 0\},\;\;i=1,\ldots,s.
\end{eqnarray}

By our assumption that the immersion is substantial, one has that
$h(M)$ is not contained in $H_1\cap\cdots\cap H_s$.

\vskip10pt

We claim that there is a sequence $(p_k)$ in $M$, $h(p_k)\notin
H_1\cap\cdots \cap H_s$, such that the distance between $h(p_k)$
and $H_1\cap\cdots \cap H_s$ tends to zero as $k\to \infty$ (the
sequence $h(p_k)$ may actually go to infinity in $\Bbb R^n$).

\vskip10pt

To prove the claim, we will need a formula for computing the
distance to the intersection $H_1\cap\dots\cap H_s$ of the affine
hyperplanes $H_1,\ldots,H_s$. Let $y$ be a fixed point in $\mathbb
R^n$. Suppose first $y\notin H_1\cap\dots\cap H_s$ and let $z$ be
the unique point in $H_1\cap\dots\cap H_s$ realizing the distance
between $y$ and $H_1\cap\dots\cap H_s$. Since $y-z\perp
H_1\cap\dots\cap H_s$, there exist unique real numbers
$a_1,...,a_s$ such that $y-z=a_1e_1+\dots+a_se_s$. Taking the
inner product with $e_j$, we obtain, for $j=1,...,s$,
\begin{eqnarray}\label{dist1}
\sum_{i=1}^sa_i\langle e_i,e_j\rangle=\langle
y-z,e_j\rangle=\langle y-p_o,e_j\rangle+\langle
p_o-z,e_j\rangle=\langle y-p_o,e_j\rangle,
\end{eqnarray}
which implies
\begin{eqnarray}\label{dist2}
a_j=\sum_{i=1}^s\langle y-p_o,e_i\rangle g^{ij},\;\;\;j=1,...,s,
\end{eqnarray}
where $(g^{ij})_{i,j=1,...,s}$ is the inverse of the matrix
$(\langle e_i,e_j\rangle)_{i,j=1,...,s}$. From (\ref{dist1}) and
(\ref{dist2}),
\begin{eqnarray}\label{dist3}
\text{d}(y,H_1\cap\dots\cap H_s)&=&||
y-z||\nonumber\\&=&\left\langle
\sum_{i=1}^sa_ie_i,\sum_{j=1}^sa_je_j\right\rangle^{\frac{1}{2}}\nonumber\\&=&\left[\sum_{j=1}^sa_j\sum_{i=1}^sa_i\langle
e_i,e_j\rangle\right]^{\frac{1}{2}}=\left[\sum_{j=1}^sa_j\langle
y-p_o,e_j\rangle\right]^{\frac{1}{2}}\nonumber\\&=&\left[\sum_{i,j=1}^s\langle
y-p_o,e_i\rangle\langle y-p_o,e_j\rangle
g^{ij}\right]^{\frac{1}{2}},
\end{eqnarray}
If $y\in H_1\cap\dots\cap H_s$, we have $\langle
y-p_o,e_i\rangle=0$, $i=1,...,s$, and (\ref{dist3}) holds in the
same way.

\vskip5pt

Assuming the claim is not true, there is $\varepsilon>0$ so that
\begin{eqnarray}\label{noclaim}
d\big(h(x),H_1\cap\cdots \cap H_s\big)\geq\varepsilon,\;\;\;x\in
M.
\end{eqnarray}
Let $H$ be the hyperplane of $\mathbb R^n$ that contains $p_o$ and
is orthogonal to the vector $e_1+\dots +e_s$, and $f:\mathbb
R^n\to\mathbb R$ the corresponding height function with respect to
$(e_1+\dots+e_s)/a$, $a=||e_1+\dots+e_s||$, so that
\begin{eqnarray}\label{heightfunction1}
f(y)=\Big\langle y-p_o,\frac{e_1+\dots+e_s}{a}\Big\rangle.
\end{eqnarray}
If $f_i(y)\geq 0,\;i=1,...,s$, and $f(y)<\delta$, it follows from
(\ref{heightfunction}) and (\ref{heightfunction1}) that
\begin{eqnarray}
0\leq \langle y-p_o,e_i\rangle <a\delta,\;\;\;i=1,...,s,
\end{eqnarray}
which implies, with the aid of (\ref{dist3}),
\begin{eqnarray}
d(y,H_1\cap\dots\cap H_s)<naC\delta,\;\;\;C^2:=\max_{i,j}|g^{ij}|.
\end{eqnarray}

Choosing $\delta=\varepsilon/naC$, we conclude that
$d(y,H_1\cap\dots\cap H_s)<\varepsilon$ for all $y\in\mathbb R^n$
satisfying $f(y)<\delta$ and $f_i(y)\geq 0,\;i=1,...,s$. It
follows from the above and (\ref{noclaim}) that
$f\big(h(x)\big)\geq\delta$, for all $x\in M$. Since the set
$\{y\in\mathbb R^n:f(y)\geq\delta\}$ is convex, we conclude that
\begin{eqnarray}
\text{Conv}\,[h(M)]\subset\{y\in\mathbb R^n:f(y)\geq\delta\},
\end{eqnarray}
contradicting the fact that $p_o$ belongs to $H$ and also to the
boundary of $\text{Conv}\,[h(M)]$. Hence (\ref{noclaim}) cannot
occur, and the claim is proved.

\vskip5pt

Since
$$
\lim _{k\to \infty} f_i\big(h(p_k)\big)= \inf _M (f_i\circ h)=0, \;\; 1\leq i \leq s,
$$
we can use Theorem  \ref{SS} to obtain $s$ sequences $q_k^{(i)}
\in M$, $1\leq i \leq s$, $k\geq 1$, such that the distance
between $q_k^{(i)}$ and $p_k$ goes to zero and $\nabla (f_i\circ
h)(q_k^{(i)})\to 0$ when $k\to \infty$. Since $\nabla (f_i\circ
h)(x)$ is the tangential component of $\nabla f_i\big(h(x)\big)$
in $T_xM$ for all $x\in M$, and $e_i=\nabla f_i(y)$ for all
$y\in\mathbb R^n$, this last condition means that the angle
between $e_i$ and the normal space $N(q_k^{(i)})$ is tending to
zero.

\vskip5pt

Passing to a subsequence, we may assume that $N(p_k)\to W$ for
some $W\in G(n-m,n)$. Since the distance between $q_k^{(i)}$ and
$p_k$ is going to zero as $k\to \infty$, and the Gauss map is
uniformly continuous, it follows that $N(q_k^{(i)})$ is also
converging to $W$.  But, as remarked before, the limit of
$N(q_k^{(i)})$ in $G(n-m,n)$ contains $e_i$. This proves that $W$
contains the $s$ linearly independent vectors $e_1,\ldots, e_s$.
In particular, $\hbox{codim} \;h(M)=\dim W =n-m\geq s$, as
desired. \qed

\vskip10pt

\noindent{\bf Proof of Theorem \ref{ConvHullC2}.} Since, by
hypothesis, the length $|\sigma|$ of the vector valued second
fundamental form $\sigma$ is bounded, the Grassmanian-valued Gauss
map is uniformly continuous, and the first assertion, $s\leq n-m$,
follows from Theorem \ref{ConvHull}.

Defining $f_i:\mathbb R^n\to\mathbb R$ by (\ref{heightfunction}),
it follows from our assumption on the vectors $e_1,\ldots,e_s$
that the functions $f_i\circ h,\;i=1,\ldots,s$, are all
nonnegative on $M^m$. A simple calculation shows, for
$i=1,\ldots,s$, that
\begin{eqnarray}\label{ConvHullC21}
\text{Hess} (f_i\circ h)_x(v,v)=\langle
\sigma(v,v),e_i\rangle,\;\;\;x\in M,\;v\in T_xM,
\end{eqnarray}
and so
\begin{eqnarray}\label{ConvHullC22}
\Delta (f_i\circ h)(x)=n\langle \overrightarrow
H(x),e_i\rangle,\;\;x\in M,\;i=1,\ldots,s.
\end{eqnarray}

Let $(p_k)$ be a sequence in $M^m$ such that $\text{d}(h(p_k),
H_1\cap\cdots \cap H_s)\to 0$ as $k\to \infty$. That such a
sequence exists is a consequence of the proof of Theorem \ref
{ConvHull}. It is immediate that $(p_k)$ is a minimizing sequence
for each one of the functions $f_i\circ h$.

Using that $|\sigma|$ is bounded, we obtain from
(\ref{ConvHullC21}) that the operator norm of $\text{Hess}
(f_i\circ h)$ is uniformly bounded on $M$, for all $i=1,\ldots,s$.
Since, by the Gauss equation, the Ricci curvature of $M$ is
bounded, we can apply Theorem \ref{weakyau} to obtain $s$
sequences $\big(q_k^{(i)}\big)$ in $M$, $1\leq i \leq s$, $k\geq
1$, such that
\begin{eqnarray}\label{ConvHullC23}
\lim_{k\to\infty}\text{d}\big(q_k^{(i)},p_k\big)=0,
\end{eqnarray}
\begin{eqnarray}\label{ConvHullC24}
\liminf_{k\to \infty}\Delta (f_i\circ h)\big(q_k^{(i)}\big)\geq 0.
\end{eqnarray}
It follows from (\ref{ConvHullC22}) and (\ref{ConvHullC24}) that
\begin{eqnarray}\label{ConvHullC25}
\liminf_{k\to \infty}\left\langle \overrightarrow
H\big(q_k^{(i)}\big),e_i\right\rangle\geq 0,\;\;\;i=1,\ldots,s.
\end{eqnarray}
Using (\ref{ConvHullC23}) and our assumption that the mean
curvature vector field $\overrightarrow H$ is uniformly continuous
on $M$, we obtain $||\overrightarrow H(p_k)-\overrightarrow
H(q_k^{(i)})||\to 0$, which implies, with the aid of
(\ref{ConvHullC25}), that $\liminf_{k\to \infty}\left\langle
\overrightarrow H(p_k),e_i\right\rangle\geq 0,\;i=1,\ldots,s$.
This concludes the proof of the theorem.\qed

\vskip10pt

\noindent{\bf Proof of Theorem \ref{ConvHullC2a}.} From the proof
of Theorem \ref{ConvHull}, there exists a sequence $(p_k)$ in
$M^m$ such that
\begin{eqnarray}\label{ConvHullC2a1}
\lim_{k\to\infty}\text{d}(h(p_k),H_1\cap\dots\cap H_s)=0.
\end{eqnarray}
Since $M^m$ is compact we may assume, passing to a subsequence,
that $(p_k)$ converges to a point $p\in M$. It follows from
(\ref{ConvHullC2a1}) that $h(p)\in H_1\cap\dots\cap H_s$, which
proves the first assertion of the theorem.

Let $q$ be a point in $M^m$ such that $h(q)\in H_1\cap\dots\cap
H_s$. Then $h(q)$ belongs to the boundary of
$\text{Conv}\,[h(M^m)]$ and $H_1,\ldots,H_s$ are supporting
hyperplanes that pass through $h(q)$. Let $(q_k)$ be a sequence in
$M^m$ so that $q_k\to q$. Then $\text{d}(h(q_k),H_1\cap\dots\cap
H_s)\to 0$, and from the proof of Theorem \ref{ConvHullC2} we
obtain
\begin{eqnarray}\label{ConvHullC2a2}
\liminf_{k\to \infty}\left\langle \overrightarrow
H(q_k),e_i\right\rangle\geq 0,\;i=1,\ldots,s.
\end{eqnarray}
Since $q_k\to q$, it follows from (\ref{ConvHullC2a2}) that
$\left\langle \overrightarrow H(q),e_i\right\rangle\geq
0,\;i=1,\ldots,s$, as was to be proved.\qed

$$
\begin{array}{lcccccccccl}
\text{Francisco Fontenele}            &&&&&&&&& & \text{Frederico Xavier}\\
\text{Departamento de Geometria}      &&&&&&&&& & \text{Department of Mathematics}\\
\text{Universidade Federal Fluminense}&&&&&&&&& & \text{University of Notre Dame}\\
\text{Niter\'oi, RJ, Brazil}          &&&&&&&&& & \text{Notre Dame, IN, USA}\\
\text{fontenele@mat.uff.br}           &&&&&&&&& & \text{fxavier@nd.edu}\\
\end{array}
$$
\end{document}